\documentclass[12pt]{article}
\usepackage{amsmath}
\usepackage{amssymb,latexsym,eucal}
\usepackage{color}
\usepackage[dvips]{graphicx}
\usepackage{authblk}

\newcommand{\ZZ}{\mathbb{Z}}
\newcommand{\NN}{\mathbb{N}}

 \newtheorem{thm}{Theorem}
\newtheorem{defns}{Definition}
\newtheorem{lem}{Lemma}

\newtheorem{prop}{Proposition}

\begin{document}


\title{Density of periodic points, invariant measures and almost equicontinuous points of  Cellular automata }

\author[1]{Pierre Tisseur }
\affil[1]{ Centro de Matematica, Computa\c{c}\~ao e Cogni\c{c}\~ao, Universidade Federal do ABC, Santo Andr\'e, S\~ao Paulo, Brasil \thanks{E-mail address: \texttt{pierre.tisseur@ufabc.edu.br}}}

\date{}

\maketitle

\begin{abstract}
Revisiting the notion of $\mu$-almost equicontinuous cellular automata introduced by R. Gilman, we show that
the sequence of image measures of  a shift ergodic measure $\mu$ by iterations of such automata  converges in
Ces\`{a}ro mean to an invariant measure $\mu_c$.
Moreover the dynamical system (cellular automaton $F$, invariant measure $\mu_c$) has  still the 
$\mu_c$-almost equicontinuity property and the set of  periodic points is 
dense in the topological support 
of the  measure $\mu_c$.
We also show that the density of periodic points in 
the topological support of a measure $\mu$
 occurs for each
 $\mu$-almost equicontinuous cellular automaton when $\mu$ is an invariant and shift ergodic measure.
 Finally using most of these results we give a non trivial  example of a
couple ($\mu$-equicontinuous cellular automaton $F$, shift  and $F$-invariant measure $\mu$) such that the 
restriction of $F$ to the topological support of $\mu$ has no equicontinuous points.
 
\end{abstract}



\section{Introduction}

A one-dimensional cellular automaton (CA) is a discrete mathematical idealization of a space-time
physical system. The space, called configuration space, is the set of
doubly infinite sequences of elements of a finite set $A$. The discrete time is
represented by the action of a CA  on this space.
Using extensive computer simulations, Wolfram in \cite{Wo86} has proposed a first empirical (visual)
classification of one dimensional CA.
In \cite{GI87} Gilman propose a formal and measurable classification by roughly dividing the
set of all CA in two parts, those with {\it almost equicontinuous
points or equicontinuous points} and those with {\it almost expansive points} (partition in order and disorder).
 The Gilman's classification is defined thanks to
  Bernoulli measures 
 and corresponds to  the Wolfram's classification based on simulations 
that use random entries.
The  measure does not need to be invariant, so the
Gilman's classification can be apply to any CA.
In \cite{Ku94}, K\accent23urka, introduce a topological classification based on the
 equicontinuity, sensitiveness and expansiveness properties.
If a CA has equicontinuous points, then there exist finite configurations that
stop the propagation of the  perturbations on the one dimensional lattice. If a CA 
has $\mu$-almost equicontinuous
 points then the probability  that a perturbation move to infinity is equal to zero (see  \cite{GI87}).
 Remark that the class of CA  with almost equicontinuous points contains the
 topological class of CA with equicontinuous points.
 In this paper, we consider the definitions of Gilman ($\mu$-expansiveness and $\mu$-equicontinuity)
in the more general case of probability measure on the configuration space $A^\ZZ$.
In this case we call $\mu$-equicontinuous CA, any CA which has a set of measure one of 
$\mu$-equicontinuous points.
Our main goal is to study  the $\mu$-equicontinuous CA when $\mu$ is an invariant measure and show the existence 
of such a measure. 
 Here we prove (see Theorem \ref{thm1})  that if  $\mu$ is a shift ergodic measure and $F$ a CA  
 which  has $\mu$-equicon\-ti\-nuous points then  the sequence $(\mu\circ F^{-n})$ converges in Ces\`{a}ro mean
 to an $F$-invariant measure $\mu_c$. 
 We also show  that this automaton $F$ is still a $\mu_c$-equicontinuous CA.
  Then,  we  describe properties of  $\mu$-equicontinuous  CA when 
  $\mu$ is an $F$-invariant measure. Here we show that (see Proposition \ref{pro4})
  under these assumptions, 
  the measure entropy is equal to zero.
  If the measure is shift ergodic or is the measure $\mu_c$ given by Theorem \ref{thm1}
    (see Proposition \ref{pro5} and  \ref{pro8}), the set of $F$-periodic points is dense in the 
    topological support of this invariant measure.     
  This result extends a previous result on the density of periodic points of surjective CA with equicontinuous points acting on a mixing subshift of finite type (see \cite{BLTI}).
 Remark that in \cite{BK} Boyle and Kitchen have shown that closing CA 
 always have a dense set of periodic points. The expansive CA and some other CA 
  with equicontinuous points  belong to this large class.   
In \cite{GI87} Gilman gives an example of a $\mu$-equicontinuous CA that has no equicontinuous points.
The invariant measure $\mu_c$, (limit by Ces\`{a}ro mean of $(\mu\circ F^{n})$)
that we can construct (using our results) for this particular automaton still has $\mu_c$-equicontinuous
 points and no equicontinuous points, but the restriction of this CA to the topological support of $\mu_c$
 always  has equicontinuous points.
Using most of our results, we describe a particular CA called  $F_e$ with a non trivial dynamic
which  keep the sensitiveness property (no equicontinuous points)
if we restrict its action to the topological support  $S(\mu_c )$ of the
invariant measure $\mu_c$.
The proofs of sensitiveness and $\mu_c$-equicontinuity of the automaton $(S(\mu_c ),F_e)$ use  
non classical arguments. 
 
\section{Definitions and preliminary results}

\subsection{Symbolics systems and  cellular automata}

Let $A$ be a finite set or alphabet. Denote by $A^*$ the set of all
concatenations of letters in $A$. These concatenations are called words. The
length of a word $u\in A^*$ is denoted by  $\vert u\vert$.
The set of bi-infinite sequences
$x=(x_i)_{i\in\ZZ}$ is denoted by $A^\ZZ$. A point $x\in A^\ZZ$ is called a
configuration. For integers $i,j$ with $i\le j$  we denote by $x(i,j)$ the word $x_i\ldots
x_j$ and by $x(i,\infty )$ the infinite sequence $(v_n)_{n\in\NN}$ such that for
all $n\in\NN$ one has $v_n=x_{i+n}$.
For each integer $t$ and each word 
$u=u_1\ldots u_{\vert u\vert}$, we
call cylinder the set $[u]_t=\{x\in A^\ZZ : x_t=u_1\ldots ;x_{t+\vert
u\vert}=u_{\vert u\vert}\}$.  
We endow $A^\ZZ$ with the product
topology of the discrete topologies on the sets $A$. 
For this topology $A^\ZZ$ is a compact metric
space. A metric compatible with this topology can be defined by the distance
$d(x,y)=2^{-i}$ where $i=\min\{\vert j\vert \,\mbox{ such that } x_j\ne y_j\}$.
The shift $\sigma \colon A^\ZZ\to A^\ZZ$ is defined by :
$\sigma (x)=(x_{i+1})_{i\in \ZZ}$. 
The dynamical system $(A^\ZZ ,\sigma )$ is called the full shift. A subshift $X$
is a closed shift-invariant subset $X$ of $A^\ZZ$.  
Consider a probability measure $\mu$ on the Borel sigma-algebra
$\mathcal{B}$ of $A^\ZZ$. If $\mu$ is $\sigma$-invariant then the topological
support of $\mu$ (which is the smallest closed subset of measure 1) is a subshift denoted by $S(\mu )$.
If $\alpha =\{A_1,\ldots ,\,
A_n\}$ and $\beta =\{B_1,\ldots ,\, B_m\}$ are two partitions of a compact space $X$, denote
by $\alpha \vee \beta$ the partition $\{A_i\cap A_j\, ; 1\le i \le n; \,\,
1\le j\le m\}$.
The metric entropy $h_\mu (T)$ of a transformation $T$
is an isomorphism invariant between two $\mu$-preserving
transformations. Put $H_\mu (\alpha ) = \sum_{A\in\alpha}\mu (A)\log \mu (A)$.
The entropy of the partition $\alpha$ is defined as $h_\mu (\alpha ) =
\lim_{n\to\infty}\frac{1}{n}H_\mu (\vee_{i=0}^{n-1}T^{-i}\alpha )$ and the entropy of $(X,T,\mu )$ as
$\sup_\alpha h_\mu (\alpha )$. 
 A cellular automaton  is a continuous self-map $F$ on
$A^\ZZ$ commuting with the shift. The Curtis-Hedlund-Lyndon theorem 
states  that for every  $F$ there exist an integer $r$ and a
block map $f$ from $A^{2r+1}$ to $A$ such that: $F(x)_i=f(x_{i-r},\ldots ,x_i
,\ldots ,x_{i+r}).$ 
 The integer $r$ is called the radius of the CA. 
 If $X$ is a subshift of  $A^\ZZ$ and one has
$F(X)\subset X$, the restriction of $F$ to $X$ determines a dynamical system
$(X,F)$; it is
called a CA on $X$.
For example, given any shift invariant measure  we can consider the restriction of 
 $F$  to $S(\mu )$.
A closed subset of $Y\subset A^\ZZ$ (not necessarily shift-invariant) such that
$F(Y)\subset Y$ is said $F$-invariant.
The omega limit set of any set $S\subset A^\ZZ$ under the action of  $F$ is denoted by 
$w(S,F)=\lim_{n\to\infty}\cap_{j=0}^n\cup_{i=j}^\infty \{F^i(S)\}$. 
 

\subsection{Almost equicontinuous points of cellular automata}

In \cite{GI87}  Gilman  shows that considering any  Bernoulli  measure $\mu$, it is possible to divide the
 set of all CA in the three following classes: The class of CA where there exist equicontinuous points,
the class of CA with $\mu$-almost equicontinuous points but without equicontinuous point and the class of
$\mu$-almost expansive CA.
In this section we recall the topological and classical  definitions for
the expansive and equicontinuous classes
of CA acting on  $A^{\ZZ}$ and we extend 
 the Gilman's measurable definitions to any probability measure $\mu$. 
 
For any integer  $n\ge 0$ and point $x\in A^\ZZ$, we denote by  $B_n(x)$ the set of  points $y$ such that for all
$i\in\NN$ one has $d(F^i(x),F^i(y))\le 2^{-(n+1)}$ and by  $C_n(x)$   
the set of  points $y$ such that $y_j=x_j$ with $-n\le j\le n$.

\begin{defns} Equicontinuity
$\mbox{ }$

-A point $x\in A^\ZZ$ is called an equicontinuous point if for all positive integer $n$ there exists another
positive integer $m$ such that $B_n(x)\supset C_m(x)$.

-A CA is called almost equi\-continuous when there exists equi\-conti\-nuous points.

-A CA is equicontinuous when all points in $A^\ZZ$ are equicontinuous.

-A point $x$ is  $\mu$-almost equicontinuous if  for all $m\in\NN$ one has 
$$
\lim_{n\to\infty}\frac{\mu \left(C_n(x)\cap B_m(x)\right )}{\mu (C_n(x))}=1.
$$

- A CA is $\mu$-almost equicontinuous if there exists a set of full measure of 
$\mu$-almost equicontinuous points.
\end{defns}

The last definition  extends the notion of $\mu$-equi\-con\-tinuity to all measurable CA $(A^\ZZ, F, \mu)$.
If $x$ is a $\mu$-equi\-con\-ti\-nuous point which belongs to the topological
support $S(\mu )$ of some  measure $\mu$ then $x$ is also a $\mu$-equicontinuous point.
Remark that for CA, to have no equicontinuous points is equivalent to sensitiveness (see \cite{Ku94}).
\begin{defns} Expansiveness
$\mbox{ }$

-A Cellular automaton is positively expansive if there exists a positive integer  $n$
such that for all $x\in A^\ZZ$ one has $B_n(x)=\{x\}$.

-A  cellular automaton  $F$ is almost expansive if there exists a positive integer  $n$
 such that for all $x\in A^\ZZ$,
$\mu (B_n(x))=0$. 
\end{defns}

In \cite{GI87} Gilman proves the following  proposition 
that allows to establish a classification based on $\mu$-equicontinuity. 
Since in \cite{GI87},  the proof of the original Proposition \ref{pro1} 
requires only the ergodicity of any  
Bernoulli measure, it  
can be extended to any shift ergodic measure.

\begin{prop}\label{pro1}\cite{GI87}
Let $\mu$ be a shift-ergodic measure and $F$ a cellular automaton with
radius $r$. The following  properties are equivalent:

\hskip  .3 true cm (i) $F$ has a $\mu$-equicontinuous point.

\hskip  .3 true cm (ii) There exists a  point $x\in A^\ZZ$ such that $\mu (B_m(x))>0$ 
for all $m\ge 0$.

\hskip .3 true cm (iii) There exists a  point $x\in A^\ZZ$ such that $\mu (B_r(x))>0$.

\hskip .3 true cm (iv) The set of $\mu$-equicontinuous points has measure 1 for $F$.
\end{prop}

When $\mu$ is a shift ergodic measure and $F$ has no $\mu$-equicontinuous points then 
$F$ is an an almost expansive CA.
A point $x$ is an equicontinuous point when the interior of $B_n(x)$ is non empty for all $n\in\NN$.

In \cite{GI87}, Gilman introduce a measure-theoretic analogue of the interior of a set.
For any measurable set $E$, define $\rho_E(x)=\lambda$ if 
$\lim_{n\to\infty}\frac{\mu (C_n(x)\cap E)}{\mu (C_n(x))}=\lambda$ and 
call  $E^{\mu}$ the set $\{x\in E:\rho_{E}(x)=1\}$.

\begin{lem}(Lebesgue)\label{lm1}
If $\mu$ is a Borel probability measure and $E$ any measurable set,  we have $\mu (E^{\mu})=\mu (E)$.
\end{lem}




A point $x$ is a $\mu$-equicontinuous point if  for all $n\in\NN$ one has $x\in B^{\mu}_n(x)$.

The next topological result is due to Gilman (see \cite{GI88}). 
We give here a detailed proof of this key result. 

\begin{prop}\label{pro2}\cite{GI88}
If there exist a point $x$ and an integer $m\neq 0$ such that $B_n(x)\cap \sigma^{-m}B_n(x)\neq\emptyset$
with $n\ge r$ (the radius of the automaton $F$) then the common
sequence {\small $(F^i(y)(-n,n))_{i\in\NN}$} of all points $y\in B_n(x)$ is ultimately periodic.
\end{prop}

{\it Proof:}\\
First remark that for each shift periodic point $\overline{w}$ of period $P$,  the cardinal of the set 
$\{F^i(\overline{w})\vert i\in\NN\}$ 
is finite and less or equal to $(\#A)^P$ (for all integer $k\ge 0$ one has 
$\sigma^{P}\circ F^k(\overline{w})=F^k\circ\sigma^{P}(\overline{w})=F^k(\overline{w})$).
This implies that the sequence $(F^{i}(\overline{w}))_{i\in\NN}$ is ultimately periodic.
Since all the elements of $B_n(x)$ share the same ultimately periodic sequence $(F^i(x)(-n,n))_{i\in\NN}$
, we only need to show  that $B_n(x)$ contains a shift periodic element $\overline{w}$ to finish the proof.
Now suppose without loosing generalities that $m>0$,  
pick a point $y_1\in B_n(x)\cap\sigma^{-m}(B_n(x))$ and put $w:=y_1(-n,-n+m-1)$.

We claim that for all points $x$, integers  $n\ge r$ and  $y\in B_n(x)$,  
all  points $z,z'$ that verify  $z(-\infty ,-n)=y(-\infty ,-n)$, 
$z(-n+1,+\infty )=x(-n+1,+\infty )$, $z'(-\infty ,-n)=x(-\infty ,-n)$ and 
 $z'(-n+1,+\infty )=y(-n+1,+\infty )$ belong to $B_n(x)$.\\ 
 We prove the claim only for $z$ using a recurrence proof. To simplify, we indifferently denote by $f$, the 
local rule of $F$ which is a block map from $A^{2r+1}$ to $A$ and also all the finite extensions of $f$ which are  block maps from 
$A^{2r+1+k}$ to $A^{k+1}$ with $k\in\NN$.
Since $y\in B_n(x)$ 
 we have $z(-n,n)=y(-n,n)=x(-n,n)$. 
Suppose that for $i>0$ one has $F^t(z)(-n,n)=F^t(x)(-n,n)$ for  $0\le t\le i$.
In this case we have $F^t(z)(-\infty , -n)$ $=F^t(y)(-\infty ,-n)$ for $0\le t\le i$ and  
 $F^{i+1}(z)(-n,0)=f(F^i(z)(-n-r,r))=f[F^i(z)(-r-n,-n-1)F^i(x)(-n,r)]=
f[F^i(y)(-r-n,-n-1)F^i(x)$ $(-n,r)]=F^{i+1}(x)(-n,0)$. 
 Since  
$F^{i+1}(z)(1,n)=f(F^i(x)(-r+1,n+r))=F^{i+1}(x)(1,n)$ it follows  that $F^{i+1}(z)(-n,n)=F^{i+1}(x)(-n,n)$. 
We can conclude saying  that $(F^k(z)(-n$ $,n))_{k\in\NN}$ $=(F^k(x)$ $(-n,n))_{k\in\NN}$ 
which implies  that $z\in B_n(x)$.

  Now we apply the claim to a point $\sigma^m (y_1)\in \sigma^m B_n(x)\cap B_n(x)$.
  The points $y_1$ and $\sigma^m (y_1)$ belong to $B_n(x)$, so the point $y_2$ 
  such that $y_2(-\infty ,-n)=\sigma^m (y_1)(-\infty$ $,-n)$ and $y_2(-n+1,+\infty )=
  y_1(-n+1,+\infty )$ belongs to $B_n(x)$.
 We can see  that $y_2(-n-m,-n+m-1)=ww$ and $y_2\in \sigma^m B_n(x)\cap B_n(x)\cap\sigma^{-m}B_n(x)$.
Next we construct $y_3$ by applying the claim to $\sigma^m(y_2)$ and $\sigma^{-m}(y_2)$ (remark that 
 $\sigma^m(y_2)\in B_n(\sigma^{-m}(y_2))=B_n(x)$). The point $y_3$ is such that 
 $y_3(-\infty , -n)=\sigma^{m}(y_2)(-\infty , -n)$ and $y_3(-n+1,+\infty )=\sigma^{-m}(y_2)(-n+1,+\infty)$. 
 Since we have $y_3(-n-2m,2m-1)=wwww$, we can repeat the same process to $\sigma^{2m}(y_3)$ and $\sigma^{-2m}(y_3)$ 
 and  construct a point $y_4$ such that $y_4(-n-4m,4m-1)=w^8$. 
 Finally, the sequence  $(y_n)_{n\in\NN}$ of points of $B_n(x)$ that we  construct  by this algorithm 
 converges  to the shift periodic point $\overline{w}=^\infty w^\infty=\ldots wwww\ldots$. 
  Since $B_n(x)$ is closed and compact, $\overline{w}$ is in $B_n(x)$.

\hfill$\Box$

In \cite{GI88} Gilman state the following result  using the ergodic properties
of any Bernoulli  measure $\mu$.

\begin{prop}\cite{GI88}\label{pro3}
Let $\mu$ be a shift ergodic measure. If a cellular automaton $F$ has a $\mu$-equicontinuous point,
 then for all  $\epsilon >0$ there exists a $F$-invariant closed set $Y$
such that $\mu (Y)>1-\epsilon$ and the restriction of $F$ to $Y$ is equicontinuous.
\end{prop}

{\it Proof:}\\
Let $x$ be a $\mu$-equicontinuous point and $m$ and integer greater than $r$, 
 the radius  of  $F$. 
 Since $\mu$ is a shift ergodic measure and $\mu (B_m(x))>0$, 
 for each  integer $k\ge 0$ and $\mu$-almost all points $y$, there exist  positive integers 
$i$, $j$ greater than $k$ such that $y\in \sigma^{-i}B_m(x)\cap\sigma^j B_m(x)$.
Then from Proposition \ref{pro2}, the sequences  
$(u_n)_{n\in\NN}=(F^n(y)(i-m,i+m))_{n\in\NN}$ and $(v_n)_{n\in\NN}=(F^n(y)(-j-m,-j+m))_{n\in\NN}$
 are ultimately periodic  with respective preperiods  $pp_u$ and $pp_v$.
For all $n\in\NN$, denote by $w_{n}$ the word $F^n(y)(-j+m+1,i-m)$ and 
remark that  for each integer $n\ge 1$, one has $w_n=f(u_{n-1}w_{n-1}v_{n-1})$ where $f$ is a block map.
  It follows that if a word $w=w_t$ with $t\ge pp_u,pp_v$ appears infinitely often in $(w_{n})_{n\in\NN}$, 
   the sequence $(w_{n+t})_{n\in\NN}$ is  periodic. 
  This implies that  $(w_n)_{n\in\NN}$ and $(F^n(y)(-k,k))_{n\in\NN}$ are  
    ultimately periodic sequences.
 Let $P$ be a map from $\NN$ to $\NN^2$ and  $Y_{P(k)}$ be the set of points $y$ such that each sequence 
 $(F^n(y)(-k,k))_{n\in\NN}$ are periodic of period $p(k)$ and preperiod $pp(k)$ where $(p(k),pp(k))=P(k)$.
As $\mu$ is shift ergodic measure, for each real $\epsilon >0$ there exists a map
$P_\epsilon\colon \NN\to\NN^2$ such that for all $n\in\NN$ we have 
$\mu (Y_{{P_{\epsilon}(n)}})>1-\epsilon\times 2^{-n}$.
Since each $Y_{{P_\epsilon}(k)}$ is closed and $F$-invariant, 
the set $Y_\epsilon =\lim_{n\to\infty}\cap_{k=1}^nY_{{P_\epsilon}(k)}$ is closed and $F$-invariant too. 
Clearly $\mu (Y_\epsilon )> 1-\epsilon$ and 
 each point $y\in Y_\epsilon$  is an equicontinuous point  since 
 for each integer $k\ge 0$  one has $C_{l\times r}(y)\subset B_k(y)$ with  
 $l=p(k)+pp(k)$ and $(p(k),pp(k))=P_\epsilon (k)$.

\hfill$\Box$
 
\section{Results on invariant measures}

\subsection{Measure entropy and density of the set of periodic points}


\begin{prop}\label{pro4}
The measure entropy $h_\mu (F)$ of a $\mu$-equicontinuous and $\mu$-invariant 
cellular automaton $F$ is equal to zero.
\end{prop}
{\it Proof:}\\
Let $\alpha_p$ be the partition of $A^\ZZ$ by the $2p+1$ central coordinates. Two points $x$ and
$y$ belong to the same element of $\alpha_p$ if and only if $x(-p,p)=y(-p,p)$.  Let
$\alpha_p^n(x)$ be the element of the partition $\alpha_p\cap F^{-1}\alpha_p\ldots F^{-n+1}\alpha_p$
 which contains $x$.
Clearly for all $n\in\NN$, we have $\alpha^n_p(x)\supset B_p(x)$.
Since almost all points are $\mu$-equicontinuous points, there exist a set $Z$ with measure 1 such that
if $y\in Z$ then $\mu (B_{m}(y))>0$ for all integer $m\ge 0$.
This implies that for almost all $y$ and all positive integer $p$, we have
$\lim_{n\to\infty}\frac{-\log \mu (\alpha_p^n(y))}{n}\le \lim_{n\to\infty}\frac{-\log \mu (B_p(x))}{n}=0$.
Using  the Shannon-Breiman-McMillan theorem which tell that
$$
h_\mu (F,\alpha_p )=\int_{A^\ZZ} \lim_{n\to\infty}\frac{-\log \mu (\alpha_p^n(y))}{n}d\mu (y)=0, 
$$
 we can conclude that $h_\mu (F)=\lim_{p\to\infty}h_\mu (F,\alpha_p)=0$.

\hfill$\Box$

In Proposition \ref{pro4}, the measure $\mu$ does not need to be shift invariant.
From Proposition \ref{pro4}, all  cellular automaton  $F$ which has  equicontinuous points
in the topological support of a shift ergodic and $F$-invariant measure $\mu$ verifies  $h_\mu (F)=0$.
This result about CA with equicontinuous points first appears  in \cite{Ti2000}.

In \cite{BRTI2007} it is shown that if a measure is shift and $F$-invariant,  
the entropy of the CA  $(A^\ZZ,F,\mu)$ is equal to zero when some discrete analogues of 
Lyapunov exponents are null.
 Remark that these Lyapunov exponents can be equal to zero for 
almost expansive CA (see \cite{BRTI2007}) but even if they are always null when there exists a set of full measure of equicontinuous points (see \cite{Ti2000}),  they are not in general equal to zero for $\mu$-equicontinuous CA 
(it can be easily seen in  examples of section 3.3).

\begin{prop}\label{pro5}
If a  cellular automaton $F$ has some $\mu$-equicontinuous points where $\mu$ is a $F$-invariant and shift ergodic 
measure then the set of $F$-periodic points is dense in the topological support of $\mu$.
\end{prop}
{\it Proof:}\\
It is sufficient to show that
 for each point $z$ in $S(\mu )$ (the topological support of $\mu$)  and positive integer $p$ we can construct a $\sigma$ and $F$ periodic point
$\overline{w}=^\infty w^\infty$ (bi-infinite repetition of the word $w$) such that $\overline{w}(-p, p)=z(-p,p)$.
 Remark that if $x$ is a $\mu$-equicontinuous point then   
   $\mu (B_r(x))>0$ where $r$ is the radius of $F$.
 Since  $z$ is in $S(\mu )$ and $\mu$ is a shift ergodic measure then $\mu (C_p(z))>0$ and 
  there exist  $(i,j)\in\NN^2$ such that
$\mu \big(C_p(z)$ $\cap\sigma^{-(i+p)}B_r(x)\cap\sigma^{j+p}B_r(x)\big )>0$.
To simplify, write $S=C_p(z)$ $\cap\sigma^{-(i+p)}B_r(x)\cap\sigma^{j+p}B_r(x)$.
Clearly there exists a point $y\in S$ such that $\mu ([y(-r-i-p,j+p-r-1)]_{-r-i-p}\cap S)>0$.
Denote by $S'$ the set $[y(-r-i-p,j+p-r-1)]_{-r-i-p}\cap S$ and 
  remark  that using the Poincar\'e recurrence theorem, we obtain that there exists an integer $m>0$ such that
$\mu\left( S'\cap F^{-m}S'\right)$ $>0$.
Remark that all points $y'\in S'$ share the same sequence
$(F^n(y')(-r-i-p,$ $j+p-r-1))_{n\in\NN}$  
since it result of the same combination of block maps  on the word $y(-i-p+r+1,j+p-r-1)$ and on the two 
sequences $(F^n(x)(-i-p-r,-i-p+r))_{n\in\NN}$ and $(F^n(x)(j+p-r,j+p+r))_{n\in\NN}$.
It follows that  $\forall y'\in S'$ one has $F^m(y')(-r-i-p,j+p-r-1)=y'(-r-i-p,j+p-r-1)$.
From  Proposition \ref{pro2} and its proof, since $\sigma^{-i}B_r(x)\cap\sigma^{j}B_r(x)\neq\emptyset$, the shift 
periodic point $\overline{w}=^\infty\hskip -.1 cm w^\infty$ such that  
$w=\overline{w}(-r-i-p,j+p-r-1)=y(-r-i-p,j+p-r-1)$ belongs to the set $S'$.
Finally, we obtain that $F^m(\overline{w})(-r-i-p,j+p-r-1)=\overline{w}(-r-i-p,j+p-r-1)$ and since  
the common $\sigma$-period of $\overline{w}$ and $F^m(\overline{w})$ is less or equal to $\vert w\vert =2p+2r+i+j$, 
we get that  $F^m(\overline{w})=\overline{w}$ 
which finish the proof.

\hfill$\Box$


\subsection{Invariant measures as limit of Ces\`{a}ro means}

Proposition \ref{pro2} also allows us to prove a Ces\`{a}ro mean convergence result.

\begin{thm}\label{thm1}
Let $\mu$ be a shift-ergodic measure. If a cellular automaton $F$ has some
$\mu$-almost equicontinuous points then the sequence $(\mu\circ F^{-n})_{n\in\NN}$ 
converges vaguely in Ces\`{a}ro mean 
to an invariant measure $\mu_c$.
\end{thm}
{\it Proof}\\
To show that the sequence of measure $(\frac{1}{n}\sum_{i=0}^{n-1}\mu\circ F^{-i})_{n\in\NN}=(\mu_n)_{n\in\NN}$
converges vaguely in measure, we need to show
that for all $x\in S(\mu)$ and $m\in\NN$ the sequence $\left(\mu_n (C_m(x))\right)_{n\in\NN}$ converges.
From Proposition \ref{pro3} there exists a set $Y_{\epsilon}$
of measure greater than $1-\epsilon$ such that for all points $y\in Y_{\epsilon}$ and positive integer $k$ 
 the sequences $(F^n(y)(-k,k))_{n\in\NN}$ are eventually periodic with preperiod $pp_\epsilon (k)$ and
period $p_\epsilon (k)$.
Hence for all $x\in A^\ZZ$ and integer $k\ge m$
$$
\begin{array}{ll}
\mu_n (C_m(x)\cap Y_{\epsilon})=&\frac{1}{n}\sum_{i=0}^{pp_\epsilon (k) -1}\mu \left
( F^{-i} \left
(C_m(x)\right )\cap Y_{\epsilon}\right )\cr
 &+\frac{1}{n}\sum_{i=pp_\epsilon (k)}^{n-1}\mu \left ( F^{-i}
\left (C_m(x)\right )\cap Y_{\epsilon}\right ).
\end{array}
$$
The first term tends to $0$; using periodicity one gets
$$
\lim_{n\to\infty}\mu_n (C_m(x)\cap Y_{\epsilon})=\frac{1}{p_\epsilon (k)}
\sum_{i=0}^{p_\epsilon (k)-1}\mu
\left ( F^{-(i+pp_\epsilon (k))} (C_m(x) \cap Y_{\epsilon }\right ).
$$

Clearly  we have
$\lim_{\epsilon\to 0}\mu_n(C_m(x)\cap Y_{\epsilon })=\mu_n(C_m(x))$.
The convergence is uniform with respect to $\epsilon$ since for all $x$ and $m\in\NN$
$$
\left\vert\mu_n (C_m(x)\cap Y_{\epsilon } )-\mu_n (C_m(x))\right\vert \le \frac{n\epsilon}{n}
=\epsilon .
$$

Consequently, letting $\epsilon$ going to $0$, we get the result by inverting the limits
$$
\lim_{n\to\infty}\frac{1}{n}\sum_{i=0}^{n-1}\mu\circ F^{-i}(C_m(x))
=\lim_{n\to\infty}\frac{1}{n}\sum_{i=0}^{n-1}\lim_{\epsilon\to 0}
\mu\circ F^{-i}(C_m(x)\cap Y_{\epsilon} )
$$
$$
=\lim_{\epsilon\to 0}\lim_{n\to\infty}\frac{1}{n}\sum_{i=0}^{n-1}
\mu\circ F^{-i}(C_m(x)\cap Y_{\epsilon } )
$$

$$
=\lim_{\epsilon\to 0}\frac{1}{p_\epsilon (k)}\sum_{i=0}^{p_\epsilon (k)-1}\mu
\left ( F^{-(i+pp_\epsilon (k))} (C_m(x)
)  \cap Y_{\epsilon }\right )=\mu_c(C_m(x)).
$$

We denote by $\mu_c$ the Ces\`{a}ro mean limit of $(\mu\circ F^{n})_{n\in\NN}$.
\hfill$\Box$

\bigskip

In the following of this section, we suppose that $\mu_c$ is the 
Probability measure which came from the Ces\`{a}ro mean of the sequence $(\mu\circ F^{-n})_{n\in\NN}$ 
when $F$ is a $\mu$-equicontinuous CA and $\mu$ is a shift ergodic measure.
 

\begin{thm}\label{thm2}
If $\mu$ is a shift ergodic measure and $F$ has  a $\mu$-equicontinuous point then $F$ is also a 
$\mu_c$-equicontinuous CA.
\end{thm}
{\it Proof}\\
We  need to show that there exists a set of measure one (for the measure $\mu_c$) 
of $\mu_c$-equicontinuous points.
 We will prove this by showing  that  
 $\mu_c$-almost all points $y$ belong $B_m^{\mu_c}(y)$ for all $m\in\NN$.
Since $F$ is a $\mu$-equicontinuous CA then there exists a point $x$
such that $x\in B^{\mu}_m(x)$ for all $m\in\NN$. 
Moreover, for all 
positive integer $k$ one has $\mu (B_k(x))>0$.

For all $m\in\NN$, define $Y_m=\cup_{i,j\in\NN^2} (\sigma^{-i-m}B_r(x)\cap \sigma^{j+m}B_r(x))$ where $r$ is 
the radius of the automaton $F$.
Since $\mu$ is a shift ergodic measure, for all $m\in\NN$ we get  
$\mu (Y_m)=1$. Then consider the omega limit sets of $Y_m$ called 
$\Omega_m:=w(Y_m,F)=\lim_{n\to\infty}\cap_{j=0}^n\cup_{i=j}^\infty F^i(Y_m)$.
Clearly for all $m\in\NN$, we have $\mu_c(\Omega_m )=1$.
Let $\Lambda (F)=w(A^\ZZ,F)$, the omega limit set of the $F$.
Since $\mu$ is a shift ergodic measure then there exists an integer $k$ such that 
$B_r(x)\cap \sigma^k B_r(x)\neq\emptyset$ and using Proposition \ref{pro2} we obtain that the sequence 
$(F^{n}(x)(-r,r))_{n\in\NN}$ 
 is ultimately periodic of period $p$.
It follows that there exist $p$ points $z_0,\ldots z_{p-1}$ such that we 
 have  $w(B_r(x),F)=\{B_r(z_l)\cap\Lambda (F) \vert \ 0\le l\le p-1\}$.
This implies that 
$\Omega_m=\cup_{z\in [z_0\ldots z_{p-1}]}
\cup_{i,j\in\NN^2} \left(\sigma^{-i-m}B_r(z)\cap\sigma^{j+m}B_r(z)\right)\cap\Lambda(F)$.
Define the set $\Omega'_m=\cup_{z\in [z_0\ldots z_{p-1}]}
\cup_{i,j\in\NN^2} \left(\sigma^{-i-m}B_r(z)\cap\sigma^{j+m}B_r(z)\right)^{\mu_c}$ $\cap\Lambda(F)$.
Since  for any measurable set $E$  one has  $\mu (E)=\mu(E^{\mu_c})$(see Lemma \ref{lm1}) and that 
we need to take off a countable number of sets of measure zero from $\Omega_m$ to obtain 
$\Omega'_{m}$, we get that $\mu_c(\Omega'_{m})=1$ for all $m\in\NN$.


Define the set $\Omega:=\cap_{m\in\NN}\Omega'_m$ and  remark  that $\mu_c(\Omega)=1$. 
Clearly,  for all $y\in \Omega$ and all $m\in\NN$, there exist $i,j\ge m$ and $z\in w(x,F)$
 such that  $y\in \left(\sigma^{-i}B_r(z)\cap\sigma^{j}B_r(z)\right)^{\mu_c}$.
Since for all $y'\in (\sigma^{-i}B_r(z)$ $\cap\sigma^{j}B_r(z))^{\mu_c}$, 
the sequence $(F^n(y')(-m,m))_{n\in\NN}$ depends only on  $y'(-m,m)$ and the common sequences  
$(F^n(y'')(-r+j,j+r))_{n\in\NN}$
and $(F^n(y'')(-r-i,-i+r))_{n\in\NN}$ where $y''$ is any point in 
$\left(\sigma^{-i}B_r(z)\cap\sigma^{j}B_r(z)\right)^{\mu_c}$,  for all $k\ge m$, we have  
$C_k(y)\cap (\sigma^{-i}B_r(z)$ $\cap\sigma^{j}B_r(z))^{\mu_c}\subset B_m(y)$.
Since $y\in C_k(y)\cap \left(\sigma^{-i}B_r(z)\cap\sigma^{j}B_r(z)\right)^{\mu_c}$, we get that  
 $y\in B_m(y)^{\mu_c}$ which finish the proof.
\hfill$\Box$


\begin{prop}\label{pro8}
If $\mu_c$ and $F$ are respectively a measure and a cellular automaton that  verify the assumptions 
of Theorem \ref{thm2} then 
 the set of $F$-periodic points is dense in $S(\mu_c )$ (the topological support of $\mu_c$).
\end{prop}


{\it Proof}\\
Let $x$ be a $\mu$-equicontinuous point.
From the proof and definitions of Theorem \ref{thm2}, there exists a finite number of  points 
$z_0,\ldots z_{p-1}\in w(\{x\},F)$ such that $\forall m\in\NN$, 
the sets $\Omega_m=\cup_{z\in [z_0\ldots z_{p-1}]}$ $
\cup_{i,j\in\NN^2} $ $(\sigma^{-i-m}B_r(z)\cap\sigma^{j+m}B_r(z))\cap\Lambda(F)$ 
have full measure with respect to the invariant measure $\mu_c$.
It follows that for each point $y\in S(\mu_c)$ and $k\in\NN$, one has 
 $\mu_c(C_k(y)\cap\Omega_m)>0$ ($\forall m\in\NN$). 
Since for each  $y\in S(\mu_c)$ and $k\in\NN$ one has $\mu_c(C_k(y)\cap\Omega_k)>0$, using the 
$\Sigma$-additivity of $\mu_c$,  we can see that 
there always exist integers $i,j\ge k$ and point $z\in w(\{x\},F)$ such that $
\mu_c ( C_k(y)\cap\sigma^{-i} B_r(z)$ $\cap \sigma^{j}B_r(z))>0$ (I). 
Finally, using the  
 final part of the proof of Proposition \ref{pro5} 
  that only use the $F$-invariance of $\mu_c$ and inequality (I), we get that there exists a dense set 
  of $F$-periodic points  in $S(\mu_c )$.

\hfill$\Box$

We remark that the measure $\mu_c$ is not necessarily shift ergodic.

\subsection{Examples of $\mu$-equi\-con\-tinuous CA
without equi\-continuous points}

In \cite{GI87} Gilman gives an example of a $\mu$-equicontinuous CA $F_s$ that has no equicontinuous points.
The automaton $F_s$ act on $\{0,1,2\}^{\ZZ}$ and is defined thank to the following  block map of radius 1.
 $$
 \begin{array}{|c|c|c|c|c|c|c|c|c|}
 *00&*01&*02&*10&*11&*12&*20&*21&*22\cr
 0&1&0&0&1&0&2&0&2
 \end{array}
 $$
 where $*$ stands for any letter in $\{0,1,2\}$.
 Considering $0$ as a background element, the 2's move straight down, 1's move to the left 
 and 1 and 2 collide annihilate each other.
In this case the measure $\mu$ is a Bernoulli measure on $\{0,1,2\}^\ZZ$ and the existence of 
$\mu$-equicontinuous points depends on the parameters $p(0),p(1),p(2)$ of this measure.
In \cite{GI87} it is shown that if $p(2)>p(1)$ then the probability that a 2 is never annihilated is positive and 
this implies that there exist $\mu$-equicontinuous points.  
Nevertheless for this interesting example, there is no Bernoulli measures $\mu$ which is 
preserving by this  automaton  $F$.
In this case a ``natural'' way to obtain an invariant measure $\mu'$ such that $F_s$ is a 
$\mu'$-equicontinuous CA is to use Ces\`{a}ro mean of image by $F_s$ of an appropriate Bernoulli measure using 
 Theorem \ref{thm1}. The dynamical system $(\{0,1,2\}^\ZZ, F_s, \mu')$ we obtain is again 
 a $\mu'$-equicontinuous CA which is sensitive (without equicontinuous points) but 
 the CA $(S(\mu'), F_s,\mu')$ where $S(\mu' )$ is the topological support of $\mu'$
 has always equicontinuous points since $S(\mu')=\{0,2\}^{\ZZ}$ and the restriction of $F_s$ to  $\{0,2\}^{\ZZ}$  
is the identity map.
 Since it has more interest to consider a dynamical  system $(S(\mu ), F, \mu)$ 
 instead of the system $(A^\ZZ, F, \mu)$ when we mix topological and measurable conditions, we will describe here 
  a CA called  $F_e$ and  a measure $\mu_c$ such that $(S(\mu_c),F_e)$ 
  is sensitive and $\mu_c$-equicontinuous.

Roughly, 
to have $\mu$-equicontinuous points but 
 no equicontinuous points  requires  that there exist some `perturbations' that can move to
 infinity but the probability that these perturbations move to infinity is equal to zero.
 One way to get theses properties for an automaton $F$ and an invariant measure $\mu_c$ 
(obtained thanks to the results of theorem \ref{thm1}), is that  $(A^\ZZ ,F)$ generates permanently `propagating 
structures'' of all type of  sizes. The ``length of life'' of the
 ``propagating structures'' depends on their sizes.
This is roughly the dynamic of the following cellular automaton $F_e$.

\subsubsection{Definition of the cellular automaton $F_e$}

The automaton $F_e$ we consider act on $X= X^1\times X^2$ where 
$X^1=\{E_0,E_1,E_2,$ $E_3,0,R,L\}^\ZZ$ and $X^2=\{0,1\}^\ZZ$.
We define $F_e$ as the composition of 3 other CA :
 $F_e=F_3\circ F_2\circ F_1$.
To simplify, we write  $\hat{E}=\{E_0,E_1,E_2,E_3\}$ and  $\overline{E}=\{0,L,R\}$.
We denote by $x=(x^1,x^2)$ any point $x\in X$  
and by  $x_i^j$  the letter in position $i$ of $x^j$ ($1\le j\le 2$).
Next we call  ${\bf 1}_{S}(x)$ the map which is equal to  one if $x\in S$ and zero
otherwise.

The automaton $F_1$ is the identity on $X^1$ and its restriction to $X^2$ came from the
following block map $f_1$ of radius 3
$$
f_1(x_{i-3}^2,x_{i-2}^2,x_{i-1}^2,x_{i}^2,x_{i+1}^2,x_{i+2}^2,x_{i+3}^2)=
{\bf 1}_{\{1\}}(x_{i-3}^2)\times{\bf 1}_{\{1\}}(x_{i-2}^2)\times{\bf 1}_{\{1\}}(x_{i-1}^2).
$$

The automaton $F_2$ is still the identity on $X^0\times X^1$ but its action  on $X^2$ depends on $X^1$.
The block map $f_2$ is defined by:
$$
f_2\left(
\begin{array}{l}
x_{i-2}^1,x_{i-1}^1, x_i^1,x_{i+1}^1, x_{i+2}^1\cr
x_{i-2}^2,x_{i-1}^2, x_i^2,x_{i+1}^2, x_{i+2}^2\cr
\end{array}
\right)
=
\left(
\begin{array}{c}
x_i^1\cr
\vee_{j=0}^2{\bf 1}_{\{E_0\}}(x_{i-j}^1)
\end{array}
\right )
$$
where $\vee_{i=0}^{2} {\bf 1}_{\{E_0\}}x_{i-j}^1$ is equal to 1 if at least one $x_{i-j}^1$ is equal to 1 and 
equal to 0 in all the other cases.

The automaton $F_3$ is the identity map on $X^2$
and  is defined thanks to a local rule $f_3$ on $X^1$.

When central coordinate $x_i$ is an element of $\overline{E}$, the block map $f_3$ is defined by the following 
rules :

$$
\begin{array}{ll}
f_3(x_{i-11}^1,\ldots x_i^1,\ldots x_{i+11}^1)&=R \mbox{ if }
x_{i-10}^1,\ldots x_i^1,\ldots x_{i+m}^1=R0^{11+m}\cr
&\mbox{ where }m=\min\{10, \min \{k-1\; | x_{i+k}\in \hat{E}\}\}\cr
&=L \mbox{ if } x_{i-m}^1,\ldots x_i^1,\ldots x_{i+10}^1=0^{11+m}L\cr
&\mbox{ where }m=\min\{10, \min \{k-1\; | x_{i-k}\in \hat{E}\}\}\cr
&=R \mbox{ if }\exists\, 0\le k,j\le 9\mbox{ such that }\cr
& x_{i-j-k}^1, \ldots , x_i^1,\ldots , x_{i+10}^1=E^*0^{k}L0^{j+10}\cr
&\mbox{ with } j+2k+1=10\mbox{ and }E^*\in\hat{E}.\cr
&=L \mbox{ if }\exists\, 0\le k,j\le 9\mbox{ such that }\cr
& x_{i-10}^1, \ldots , x_i^1,\ldots , x_{i+j+k}^1=0^{10+j}R0^{k}E^*\cr
&\mbox{ with } j+2k+1=10 \mbox{ and }E^*\in\hat{E} .\cr
\end{array}
$$
If the central coordinate $x_i$ is an element of $\hat{E}$, for each  $i\in\{0,1,2,3\}$, one has 
$$
\begin{array}{ll}
f_3(x_{i-10}^1,\ldots ,E_{i},\ldots x_{i+10}^1)&=E_{i+1} \mbox{ if } x_{i-k}^1, \ldots, x_i^1=R0^{k-1}E_i
\mbox{ with } 0\le k\le 9 \cr
&\mbox{  }\mbox{ }\mbox{ }\mbox{ where the addition `$i+1$' is  modulo 4.}
\end{array}
$$
For all the other cases where the central coordinate $x_i$ is an element $E^*$ in
$\hat{E}$, 
we have  $f_3(x_{i-10}^1,\ldots ,E^*,\ldots x_{i+10}^1)=E^*$.
In all  other cases  not  described above we have, 
$f_3(x_{i-10}^1,\ldots $ $,x_i^1,\ldots x_{i+10}^1)=0$.

\subsubsection{The invariant measure $\mu_c$}

Let $S^1$ be a subshift of finite type  
defined by the following list of accepted words:
 words of type $LE_j0^k$ with $(1\le j\le 3 \mbox{ and } k=170)$; 
 words of type $0E_l0^k$ with $(0\le l\le 3 \mbox{ and } k=170)$;
 words of type $0^l$ with $l\in\NN$. 
A typical configuration in $S^1$ is
 
$\ldots 00E_{2}000...00LE_{0}00...00LE_{1}00...00E_{1}00\ldots 00E_300... 00LE_{0}00...$

We denote by $\mu_1$ the Parry measure on  $S^1$ (see for example \cite{CP}),  
 by $\delta_{^\infty 0^\infty}$  the Dirac measure on the point  
 $^\infty 0^\infty=\ldots 00000\ldots$ and 
 we call $\mu_I$ the measure  $\mu_1\times \delta_{^\infty 0^\infty}$ on $X$.
Since $\mu_I$ is a shift ergodic measure, using Proposition \ref{pro1}, Theorem \ref{thm1} and 
\ref{thm2}, $F_e$ is a $\mu_c$-equicontinuous CA 
 (with $\mu_c=\lim_{n\to\infty}\frac{1}{n}\sum_{i=0}^{n-1}\mu_I\circ F_e^{i}$  ) if  
there exists a point $x$ such that $\mu_I(B_r(x))>0$ .  
Remark that $S^1\times ^\infty 0^\infty$ is not an invariant set for  $F_e$.
In Section 3.3.3 we will characterize the topological support of $\mu_c$ by describing the action 
of $F_3$ on the non-invariant set $S^1$.

\subsubsection{The dynamic of $F_e$}
In this subsection, we describe the global dynamic of $F_e$ by showing in first place the
contribution of each of the 3 cellular automata $F_1$, $F_2$ and $F_3$.

{\sc The dynamic of $F_1$}

The action of $F_1$ on $X^2$ is only the shift of consecutive sequences   of letters ``1'' 
 (that we call ``trains of 1'') of one
coordinate to the right and the destruction of the two last letters ``1'' at the left side of this {\it train}.

Action of $F_1$ on a  configuration of $X^2$:
$$
...01111110000... \hskip .1 cm\raisebox{.4 cm}{$F_1$}\hskip -.55 cm\mapsto
...00001111000... \hskip .1 cm\raisebox{.4 cm}{$F_1$}\hskip -.55 cm\mapsto ...00000001100... 
\hskip .1 cm\raisebox{.4 cm}{$F_1$}\hskip -.55 cm\mapsto ...0000000000... 
$$
Remark that a {\it train }of 1 with a  length $2k+1$  will move of $k$ coordinates to
the right  before collapsing.

\medskip

{\sc The dynamic of $F_2$}

The action of the cellular automaton $F_2$ is to 'create'  a sequence of three letters `1' in
$X^2$ when there is a letter $E_0$ in $x_i^1$.

Example:
$$
x=
\begin{array}{lll}
x^1\hskip .2 cm\ldots E_0**\ldots &\hskip .1 cm\raisebox{-.4 cm}{$F_2(x)$}\hskip -.6 cm
\raisebox{-.75cm}{$\mapsto$}\hskip .5 cm \ldots E_0**\ldots &
\hskip .1 cm\raisebox{-.4 cm}{$F_2(x)$}\hskip -.6 cm\raisebox{-.7 cm}{$\mapsto$}\hskip .5 cm\ldots E_0**\ldots\cr
 x^2\hskip .2 cm\ldots \;0\;\;0\;0\;0\;0\ldots &\hskip 1.6 cm \ldots \;1\;1\;1\;0\;0\ldots&\hskip 1.6 cm 
 \ldots \;\;1\;1\;1\;0\;0\ldots  
\end{array}
$$
The symbol * replace any letter in $\{0,L,R\}$.

{\sc Action of $F_2\circ F_1$}

If $F_e^i(x^1)=E_0$ for $0\le i\le n-1$, then there is at least a  {\it train of 1} of length $n+3$ moving
to the left of at least  $\lceil \frac{n+3}{2}\rceil$ coordinates.

{\sc Example of the action of $F_2\circ F_1$ (on $x=\raisebox{-.16 cm}{$x^1$}\hskip -.4 cm\raisebox{.16 cm}{$x^2$}$\hskip .2 cm  
} where $*$ replace any letter in $\hat{E}\cup \overline{E}$)
$$
\begin{array}{lllr}
x^1=^\infty 0E_0*******\ldots  &\raisebox{-.5 cm}{$(F_2\circ F_1)^n$}\hskip -1.4 cm
\raisebox{-.85cm}{$\Longrightarrow$}
&\hskip .6 cm^\infty 0E_0*******\ldots &=F^n (x)^1\cr
x^2=^\infty 0\;\;0\;0\;0\;0\;0\;0000\cdots & &\hskip .6 cm 
\underbrace{^\infty 0\;\;1\;\;1\;\;1\;\;1\cdots 1}_{n+3 \mbox{ times}}000\cdots  &=F^n(x)^2
\end{array}
$$


{\sc The oscillator dynamic of $F_3$ and the subshifts $S^1$ and $S(\mu_c)$ }

First remark that $S(\mu_c)\subset w(S^1\times ^\infty 0^\infty ,F_e)$ where 
$w(S^1\times ^\infty 0^\infty ,F_e)=\lim_{n\to\infty}$ 
$\cap_{j=0}^n\cup_{i=j}^\infty F_e^i(S^1\times ^\infty 0^\infty)$.\\
 When we apply $F_3$ to a finite configuration of the subshift $S^1$,
it is easy to see that only finite configurations of the type $0^k LE_j0^l$ will be affected by the first iteration 
of the automaton $F_3$ where $k,l\in\NN^2$ and $j\in [0..3]$. 
We have $F_3(\ldots 0^k LE_j0^{l+m}\ldots)=\ldots L0^{10}E_j0^l\ldots$ and 
$F_3^2(\ldots 0^k LE_j0^l\ldots)=\ldots L0^{20}E_j0^l\ldots$ where  $m,k\ge 170$, $l\in\NN$ and $0\le j\le 3$.
Remark that the $L$ moves to the left indefinitely unless it would  approach a letter $E^*$ in $\hat{E}$.
If there is less than nine  $0$ between the $L$ and the $E^*$, the $L$ disappears and a $R$ appears 
 at $i$ coordinates to the right of $E^*$ if the $L$ was at $10-i$ coordinates to the right of the $E^*$. 
An similar  process occurs when after some iterations, the $R's$ appears and move to the right until they encounter 
 a letter $E^*\in\hat{E}$.
The letter $R$ disappears, it appears a $L$ which return to the left and the letter 
$E^*=E_j$ in the neighborhood become 
 $E_{j+1 \mbox{ mod 4}}$. 
 
Let's see a typical  evolution of an {\it oscillator}.

\small 
$$
\hskip -2.5 cm
\ldots\underbrace{E^*0\ldots 0R0\hskip -.55 cm\raisebox{.6 cm}{...$\hookrightarrow$}\ldots 0E_i...}_
{\mbox{oscillator size $50n$}}\hskip .2 cm
\raisebox{.5cm}{$F_3^{n}$}\hskip -.7cm \Rightarrow\hskip -.1 cm \ldots\underbrace{E^*0\ldots R0^kE_i}_
{\mbox{ same oscillator}}\ldots \hskip .2 cm
\raisebox{.5cm}{$F_3$}\hskip -.7cm \Rightarrow\hskip -.1 cm 
\ldots\underbrace{E^*0\ldots 0L0^lE_{i+1\mbox{ mod 4}}}_{\mbox{ same oscillator}}
\ldots\hskip -3.3 cm\raisebox{.6 cm}{$\hookleftarrow ...$}
$$
\normalsize
where  $k\le 9$, $n\ge 4$ and $l+k=10$. 
        
 We will describe only the projection on $X^1$ of the subshift 
 $S(\mu_c)=w(S^1,F_e)\times w(^\infty 0^\infty ,F_e)$.
  In $w(S^1,F_e)$: (i)
  there is at least 170 letters in $\overline{E}$ between two $E^*$.
   (ii) There is at most one letter 
 $M\in \{L,R\}$ between two $E^*$.
 (iii) At the right side of the last  $E^*$ to the right (if it exists ) there are only letters 0.
  (iv) (iii) At the left side of the last  $E^*$  (if it exists ), there are only letters 0.
 (v) There is no configurations of the type $E^*0^mE_0*^l E^*$ where $*$ is in place 
 of any letter in $\{0,R,L\}$ and $l\ge 170$. 
 We remark that the configurations of $S(\mu_c )$ can be generated by a finite automata  
 which means that $S(\mu_c )$  is a sofic  subshift.
   

{\sc Action of $F_e$ } 

Since the action of  $F_e$ on letters  $E^*\in \hat{E}$ is the identity or a permutation in $\hat{E}$,   
 the set of configurations that contains an infinite number of letters 
 in $\hat{E}$ has measure one in $S(\mu_c)$.
 
Under the action of $F_2\circ F_1$, a configuration of the type 
$\ldots E_0\ldots\times ^\infty 0^\infty$
 will generate in $X^2$ (at the same coordinate that $E_0$) a {\it train of 1} 
until the $E_0$ will change in $E_i$ (with $i\ge 1$).
Under the action of $F_e$ (see action of $F_2\circ F_1$) 
a configuration $\ldots E^*0^lM0^k E_0\ldots\times ^\infty 0^\infty$ ($E^*\in\hat{E}$, $M\in\{R,L\}$) will produce 
 a {\it train of 1} of length $\mathcal{L}$ such that $\lfloor\frac{m}{5}\rfloor +2\le \mathcal{L}\le \lfloor\frac{m}{5}\rfloor +4$ 
 with $m=l+k+1$.  In the following we will choose $\lfloor\frac{m}{5}\rfloor +2$ or $\lfloor\frac{m}{5}\rfloor +4$ 
for the length $\mathcal{L}$ of this  {\it train} of 1 according to the context. 

Let see a typical dynamic of  an {\it oscillator transmitter}

\small 
$
$\mbox{ }$\hskip 0 cm 
\raisebox{-.5 cm}{x= }
^\infty 0\ldots \underbrace{E^*0\ldots \raisebox{.3 cm}{$\hookrightarrow$}\hskip -.5 cm \ldots 
R0^9}_{\mbox{ size l}}
E_3\ldots \ldots \hskip .6 cm  \raisebox{-.3 cm}{$F_e$}
\hskip -.6 cm\raisebox{-.6 cm}{$\Longrightarrow$}
 \hskip .8 cm ^\infty 0\ldots \underbrace{E^*0\ldots \raisebox{.3 cm}{$\hookleftarrow$}\hskip -.5 cm \ldots
  L0E_0}_{\mbox{size l}}
 \ldots\ldots \hskip .6 cm\\
$\mbox{ }$\hskip .65 cm  ^\infty 0\ldots 0000\ldots\ldots \ldots  00000\ldots \hskip 1.85 cm 
 ^\infty 0\ldots 0000\ldots\ldots \ldots  00000\ldots 
$
\\
\\
$
$\mbox{ }$\hskip 0 cm 
\raisebox{-.5 cm}{$=$ }
^\infty 0\ldots \underbrace{E^*0 \raisebox{.3 cm}{$\hookleftarrow$}\hskip -.5 cm \ldots L0}_{\mbox{ size l}}
E_0\ldots \ldots \hskip .6 cm  \raisebox{-.3 cm}{$F_e^{\lceil (l/5\rceil}$}
\hskip -1 cm\raisebox{-.6 cm}{$\Longrightarrow$}
 \hskip .8 cm ^\infty 0\ldots \underbrace{E^*0\ldots \raisebox{.3 cm}{$\hookleftarrow$}\hskip -.5 cm \ldots
  L00E_1}_{\mbox{size l}}
 \ldots\ldots\ldots\ldots \hskip .5 cm\\
$\mbox{ }$\hskip .4 cm ^\infty 0\ldots 0000\ldots \ldots  0000\ldots \hskip 2.15 cm 
^\infty 0\ldots 000\ldots\ldots  0000\underbrace{111\ldots 11}_{\lfloor\frac{l}{5}\rfloor +4}00
\hskip -1.2 cm\raisebox{.4 cm}{$\cdots\Rightarrow$}\hskip .4 cm\ldots
$
\normalsize

We  call {\it oscillator transmitter} of size $l+m+1$, any pattern of the form $E^*0^lM0^mE^*$ and
 {\it void oscillator} any pattern  of the form $E^*0^k E'^*$ where $k,l+m+1\ge 170$, each $E^*$ belongs to 
 $\hat{E}$, each  $E'^*\in\{E_1,E_2,E_3\}$ 
 and $M\in\{R,L\}$. Remark that {\it oscillators} of the type $E^*0^k E_0$ does not belong to the language 
 of $S(\mu_c)$.

 For each $l\in\NN$, denote by $\mathcal{C}_l[i]$  the union of all the sets
 $[U]_i\times X^2\subset X$
where $U=[E^{*}0^j M0^k E^{*}]_i$ is a cylinder in $X^1$, $M$  replaces one letter in $\{L,R\}$,
each $E^{*}$ are any letters in $\hat{E}$ and $j,k,l$ verify  $j+k+1=l$.
Let  $\overline{\mathcal{C}_k[i]}$ be the union of  sets
 ($[E^*0^k E'^*]_{i}\times X^2$) where $E^*$ replace any letter in $\hat{E}$, $k\ge 170$ and 
 $E'^*\in\{E_1,E_2,E_3\}$ .
We call respectively {\it oscillators transmitter } in position $i$ and {\it void oscillators} in positions $i$ 
the sets
$\mathcal{C}_l[i]$ and $\overline{\mathcal{C}_l[i]}$.
Remark that $\mu_I(\mathcal{C}_l[i])=\mu_I([E^*0^{i-1}RE^*]\times ^\infty 0^\infty)$.   
Without taking into consideration the position of the {\it oscillators}, we will call respectively $\mathcal{C}_{l}$ and 
$\overline{\mathcal{C}_{l}}$ the {\it oscillators transmitter} and {\it void oscillators} of size $l$.

\begin{figure}[!h]
\centering
\mbox{ 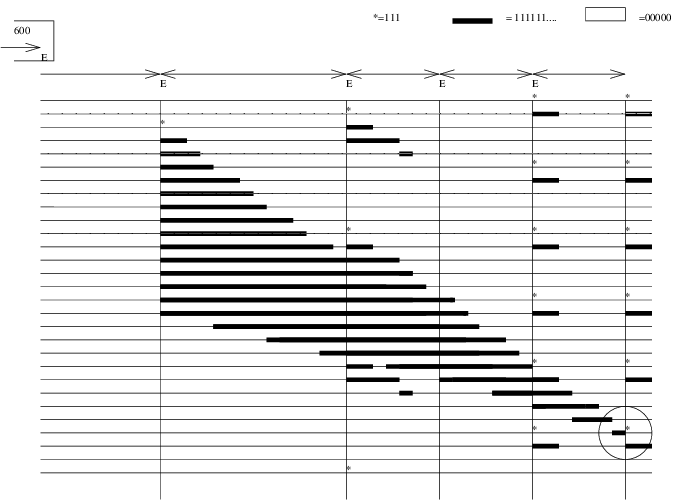 }
\caption{
An illustration of the dynamic of $F_e$ on 5 {\it oscillators} 
$\mathcal{C}_{4200}$, $\mathcal{C}_{600}$ and three $\mathcal{C}_{300}$ and the resulting dynamic of
{\it train of 1} in $X^2$.
For simplification, we do not specify the states of the {\it oscillators} and their evolution because
the interesting  part of their dynamic  can be deduce from the evolution of the {\it train of 1} in $X^2$. 
Each line represents the sequence of images of $x^2$ after every $60$ iterations of $F_e$. 
 The black horizontal lines represent the
  {\it trains of 1} and blank ones, the sequences of $0$.
  The extremity of the {\it oscillators} are delimited by arrows.
 The first {\it oscillator} of size 300 is a {\it void oscillator} which never generates any {\it trains} of 1
in $X^2$. To the left side, there is one large {\it oscillator transmitter} of size 6000 which is in a non emitting
state (the last letter is not an $E_0$).
 Remark that the non `emitting period' of the {\it oscillators transmitter  } last 3 times more than 
 the `emitting' one.
 The circle shows the end of the propagation of the train generated by $\mathcal{C}_{4200}$. 
  }
\label{fig1}
\end{figure}

{\sc Propagation of trains of 1 generated by oscillator transmitter} 


Since an {\it oscillator transmitter} $\mathcal{C}_l$ generate a {\it train of 1} with a 
length $\mathcal{L}$ at most equal to $\lfloor\frac{l}{5}\rfloor +4$ and that {\it train}  looses 2 elements when it  
moves of 1 coordinate,  it can influences some patterns  situated at
 $(\lfloor\frac{l}{5}\rfloor +4)+(\lfloor\frac{l}{5}\rfloor +4)/2=
\lfloor\frac{3l}{10}\rfloor +6$ coordinates to the right of the right extremity of the {\it oscillator} if there is no concatenation process with another train created by other  {\it oscillators transmitter}.
 Since the proofs of Proposition \ref{pro10} and \ref{noneq} only require 
 the understanding of the  propagation of  {\it trains of 1} when the initial configuration is in 
 $S^1\times ^\infty 0^\infty$,  we will only consider concatenation process with trains generated by  
 {\it oscillators transmitter} situated to the right side.  
  For two consecutive oscillators $\mathcal{C}_l\mathcal{C}_m$,  
 a  {\it trains of 1} generated by the first one  $\mathcal{C}_l$ will  
 reach the coordinates situated under the beginning of $\mathcal{C}_m$ only if 
 $m\le \lfloor\frac{3l}{10}\rfloor +6$ (iq1).
  When a {\it train } will cross the $l$ 
coordinates of another {\it oscillator}, it will loose $2l$ elements and can incorporated at most 
$\lfloor 2l/5\rfloor +8$ others ($\lfloor l/5\rfloor +4$ at the beginning and $\lfloor l/5\rfloor +4$ the end of 
the {\it train}).     
 Now consider a sequence of three {\it oscillators} $\mathcal{C}_l\mathcal{C}_m\mathcal{C}_n$ 
 where $m$ and $n$ are fixed  and $l$ is the minimum size  of the first 
 {\it oscillator} in order that its train will reach $\mathcal{C}_n$.
 The {\it train } of length at most $\lfloor l/5\rfloor +4$ will gain at most $\lfloor 2m/5\rfloor +8$ and will loose at least  $2m+2(n-(\lfloor m/5\rfloor +4))$  elements when the head of the train has  
 crossed the two {\it oscillators} (we take into consideration that  the head of the train may eventually 
  progress of $\lfloor m/5\rfloor +4$ coordinates when it passes under $\mathcal{C}_m$). 
We obtain that $l\ge 6m+10n-100$.
More generally for all sequences of $n+1$ consecutive {\it oscillators} 
$\mathcal{C}_{l_n}\mathcal{C}_{l_{n-1}}\ldots\mathcal{C}_{l_0}$, any change of the state 
of the first {\it oscillator} $\mathcal{C}_{l_n}$ will affect the {\it train of 1} situated in the 
{\it oscillator} $\mathcal{C}_{l_0}$ if  $l_n\ge \sum_{i=n-1}^{1}6l_i+10l_0-80n+60$ (iq2).

Recall that 170 is the minimum size of an {\it oscillator}.
This minimum size is required to simplify the proof of Proposition \ref{noneq} which uses 
quantitative arguments on ''flows
of trains of 1''.


\subsubsection{The topological and measurable properties of $F_e$}

\begin{prop}\label{pro10}
The dynamical system $(S(\mu_c ),F_e)$ is a  $\mu_c$-equicontinuous cellular automaton.
\end{prop}
{\it Proof}\\
  From the discussion of section 3.3.2,  we only need to show that there exists a point $x$ and an integer 
$m\ge r$  such that $\mu_I(B_m(x))>0$.
Remark that since $\mu_I(X)=\mu_I(S^1\times ^\infty 0^\infty )$  
 we will take into consideration only configuration in $S^1\times ^\infty 0^\infty$ in this proof. 
 Let $x_0=(^\infty 0^\infty,^\infty 0^\infty)$  
 and for each integer $k\ge 1$, pick
a point $x_k\in\overline{\mathcal{C}_k[-k-1-r]}$ $\cap B_r(x_0)\cap X^1\times ^\infty 0^\infty$.
We will prove that there exist integers $k>0$ such that
$\mu_I (B_r(x_k))>0$ by showing that
 $\mu_I(\overline{\mathcal{C}_k[-k-1-r]}\cap B_r(x_k)^\complement)<\mu_I (\overline{\mathcal{C}_k[-k-1-r]})$ 
 (where $B_r(x_k)^\complement$ is the complement of $B_r(x_k)$). 
The set $\overline{\mathcal{C}_k[-k-1-r]}\cap B_r(x_k)^\complement$ is the set of points  that contains
   {\it oscillators transmitter} in the left side of $\overline{\mathcal{C}_k[-k-1]}$ that are  able to 
   generate {\it trains of 1} which
move to the right and cross this {\it void oscillator} $\overline{\mathcal{C}^*_k[-k-1-r]}$ 
(the {\it trains} of ``1''  may enter in the central
coordinates ($[-r,r]$) and in this case the point does not belong to $B_m(x_0)=B_m(x_k)$).
Now, consider an {\it oscillator transmitter} of size $l$ in position $-p-l-k-r$
: $\mathcal{C}_l[-p-l-k-r]$.
Denote by ${\bf S}(p,k)$ be the minimum size of the {\it oscillator} $\mathcal{C}_l[-p-l-k-r]$ (whose
the right extremity  is situated at $p$ coordinates to the left of $\overline{\mathcal{C}_k[-k-1]}$)
in order that it
 can produce {\it trains of 1} that cross completely the {\it void oscillator} 
 $\overline{\mathcal{C}_k[-k-1]}$.
From the discussion in section 3.3.3 about the propagation of {\it train} of 1,  
 it follows that $\overline{\mathcal{C}_k[-k-1-r]}\cap B_r(x_k)^\complement$
is equal to 
 $\mathbb{S}_k=\left\{\cup_{i={\bf S}(0,k)}^\infty\mathcal{C}_i[-i-k-1-r]
\cup_{p=l_{\bf m}}^\infty\{\cup_{j={\it S}(p,k)}^\infty \mathcal{C}_j[-j-p-k-1-r]\}\right\}
$ $\cap\overline{\mathcal{C}_k[-k-1-r]}$ where 
$l_{\bf m}=170$ be the minimum size of the {\it oscillators} in $S(\mu_I)$.

 Next we  claim  that there exists a real $M\ge 0$ such that 
 for all integers $k\ge l_{\bf m}$ one has  
$
\mu_I(\overline{\mathcal{C}_k[-k-1-r]}\cap B_r(x_k)^\complement )=\mu_I(\mathbb{S}_k)
\le \mu_I (\overline{\mathcal{C}_k[-k-1-r]})\times\eta(k) \times M
$ with $\lim_{k\to\infty}\eta (k)=0$. 
Remark that since $\mu_I (\overline{\mathcal{C}_k[-k-1-r]})$ $>0$, the proof of this claim will finish the proof.

From section 3.3.3  a {\it train} of 1 generated by an {\it oscillator } $\mathcal{C}_{l_n}$ will 
cross the $n-1$ {\it oscillator } $\mathcal{C}_{l_{n-1}}\ldots\mathcal{C}_{l_1}$ and reach 
$\overline{\mathcal{C}_{l_0}}$ if $l_n\ge \sum_{i=n-1}^{1}6l_i+10l_0-80n+60$ (iq2).
 Remark that  the  last {\it oscillator} is $\overline{\mathcal{C}_k[-k-1-r]}$ which implies $l_0=k$ and 
 that the number of {\it oscillators} $n\le \frac{p}{l_{\bf m}}=\frac{p}{170}$.  
  Using (iq2)  we obtain that 
  ${\bf S}(p,k)\ge 6p+10k-80(\frac{p}{l_{\bf m}})+60\ge 5p+10k+60$ (iq3).

Since $\mu_I$ is the product of the Parry measure on the mixing subshift of finite type $S^1$ (see \cite{CP}) 
and the Dirac measure on 
$^\infty 0^\infty$, 
 then there exist a real $0<q<1$ and a 
positive integer $\mathcal{L}$ such that 
$\forall m\ge 0\;\mbox{ one has } \mu_I([u_0\ldots u_m]\times ^\infty 0^\infty)\le 
q^{\lfloor\frac{m+1}{\mathcal{L}}\rfloor}$.
   
To prove the claim put  
$\eta (k)=\mu_I\left(\mathcal{C}_{{\bf S}(0,k)}[-S(0)-k-1-r]\right)$ and using the lower bound (iq3) 
for ${\bf S}(0,k)$, we obtain  that $\eta (k)\le\mu_I\big(\mathcal{C}_{10k+60}[-10k-60-$ $k-1-r]\big)
\le q^{\lfloor\frac{11k+59}{\mathcal{L}}\rfloor}$ 
 which implies that  $\lim_{k\to\infty}\eta (k)=0$.
  Using again (iq3) for ${\bf S}(p,k)$, we obtain  that
$\mu_I(\cup_{i={\bf S}(0,k)}^\infty\mathcal{C}_i[-i-k-1-r])\le \eta (k)\sum_{i=1}^{+\infty}
q^{\lfloor\frac{i}{\mathcal{L}}\rfloor}$.
Remark that  the last and the following inequality will give us an upper bound  for the measure of $\mathbb{S}_k$:  
$$
\mu_I\left(\cup_{p=l_{\bf m}}^\infty\{\cup_{j={\bf S}(p,k)}^\infty \mathcal{C}_j[-j-p-k-1-r]\}\right)
\le  \eta (k)\sum_{i=1}^{+\infty}q^{\lfloor\frac{i}{\mathcal{L}}\rfloor}
(\sum_{j=1}^{+\infty}q^{\lfloor\frac{j}{\mathcal{L}}\rfloor}). 
$$
It follows that  we can prove the claim and consequently  finish the proof taking 
$$
M=\sum_{i=1}^{+\infty}q^{\lfloor\frac{i}{\mathcal{L}}\rfloor}+ 
\sum_{i=1}^{+\infty}q^{\lfloor\frac{i}{\mathcal{L}}\rfloor}
(\sum_{j=1}^{+\infty}q^{\lfloor\frac{j}{\mathcal{L}}\rfloor})< +\infty.
$$
\hfill$\Box$

\begin{prop}\label{noneq}
The cellular automaton $(S(\mu_c ),F_e)$ is sensitive (has no equi\-conti\-nuous points in the
topological support $S(\mu_c )$).
\end{prop}

{\it Proof}\\
If we suppose that there exists an equicontinuous point $x\in S(\mu_c)$, 
 there must exist an integer $m$ such that
$C_m(x)\subset B_0(x)$.
First suppose that if for all positive integers $n$, there exist integers $i>n$,  such that
 $(F^i(x))^2_0=0$, then  there exists $y\in C_m(x)$ such that
 $(F^i(y))^2_0=1$ (it is always possible to choose a point $y\in C_m(x)$ that contains a 
 large enough {\it oscillator transmitter} at the left side of $-m$ that  sent a {\it train of 1} that can arrived at   
coordinate $0$ after $i$ iterations).  
  This contradict the hypothesis $C_m(x)\subset B_0(x)$ and 
 it follows  that if there exists $x\in X$ and
 $m>0$ such that  $C_m(x)\subset B_0(x)$, then  there exists a positive integer 
 $n$ such that for all $y\in C_m(x)$ and $i\in\NN$
one has  $\left(F^{i+n}(y)\right)^2_0=1$  (condition (*)).
\hskip .1 cm 
Then suppose that there exists $m\in\NN$ such that $C_m(x)\subset B_0(x)$.
Since the point $z:=^\infty\hskip -.16 cm 0^\infty (-\infty ,-m-1)x(-m,\infty )$ belongs to $C_m(x)$ 
and by hypothesis   $(F^{i+n}(z))^2_0=1$ for all $i\in\NN$, 
 we obtain that  the pattern 
$x(-m,-1)$ contains  a finite sequence of  consecutive 
 {\it oscillators } $\mathcal{C}_{l_k}, \mathcal{C}_{l_{k-1}}, \ldots \mathcal{C}_{l_0}$ that 
 generate a '' continuous  flow of letters 1''. To finish the proof we need to show that this finite sequence of 
 {\it oscillators } does not exists. In order to do that, we will consider the propagation of {\it trains of 1} 
generated by a finite sequence of $k$ {\it oscillators } ($k\in\NN$).

\begin{figure}[!h]
\centering
\includegraphics[width=0.8\textwidth,height=0.26\textheight]{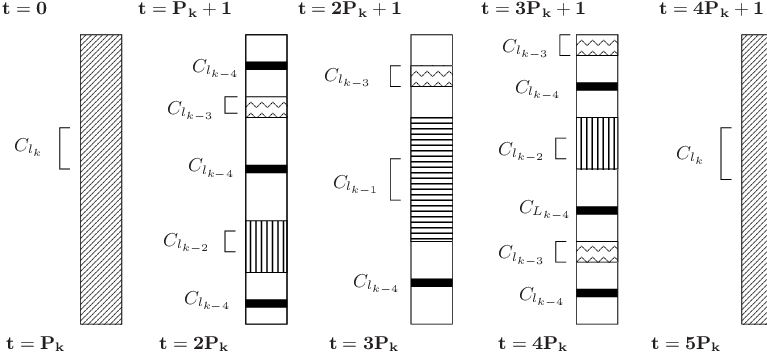}
\caption{The flow of 1 generated by a sequence of $k$ {\it oscillators}. 
The  parts not in white color represents the successive flows of letters 1 due to the trains which come from  
 the  oscillators transmitters $\mathcal{C}_{l_k}, \mathcal{C}_{l_{k-1}}, \ldots \mathcal{C}_{l_{k-4}}$ .}
\label{fig2}
\end{figure}

Let $P_i$ $(0\le i\le k)$ be an upper bound  of the number of iterations needed  for a {\it train of 1}  
 generated by the 
{\it oscillator } $\mathcal{C}_{l_i}$ to cross completely the coordinate $0$. We consider that there  is no 
concatenations with other {\it  trains} coming from the left (but possible concatenations with {\it trains} 
generated by {\it oscillators } situated to the right. 
The value of $P_i$ depends on the difference between the lost and gained letters when the left extremity of the {\it train} 
arrived in coordinate 0. It follows that $P_i\ge\lfloor\frac{l_i}{5}\rfloor +4-s_i$ where  
$s_i=\sum_{j=0}^{i-1}2l_j-\sum_{j=0}^{i-1}2(\lfloor\frac{l_j}{5}\rfloor +4)=\sum_{j=0}^{i-1}\lfloor\frac{8}{5}l_j\rfloor -8i$.
Remark that $s_i\ge 0$ if $\forall 0\le j\le i-1$, $l_j\ge 5$.  
As $l_{\bf m}=170$ we obtain that $P_i\le \lfloor\frac{l_i}{5}\rfloor +4$ for all $0\le i\le k$.
 Without loosing generalities, we can suppose that the {\it train of 1} generated by the first {\it oscillator}
 of size $l_k$
arrive in coordinate 0 at $t=0$ and last at most $P_k$ iterations. For time $t=P_k+1$ to
$4\times P_k$, there is no {\it train of 1} due to this first {\it oscillator} that pass through the
central coordinate (because $P_i\ge \lfloor\frac{l_k}{5}\rfloor +2$ and the {\it oscillator} $\mathcal{C}_{l_k}$ 
stop to generate {\it train of 1} for at least $3(\lfloor\frac{l_k}{5}\rfloor +2)$ iterations).
 The {\it train} of 1 generated by the second oscillator from the left : $\mathcal{C}_{l_{k-1}}$ last at
most $P_{k-1}$ iterations and its effect stops for a period of $3P_{k-1}$ in the interval
time $t=P_k+1$ to  $t=4P_k$.
Clearly, if $P_{k-1}$ is small enough,  between $t=P_k+1$ and  $t=4P_k$, there is  at least one interval
of length at least $3P_{k-1}$  that will be not affected by the two
first {\it oscillators}.  This interval is minimum when the {\it train }
of 1 generated by the second {\it oscillator} pass exactly in the middle of the interval
$[P_{k}+1, 4P_k]$, between 2 {\it trains} of the first {\it oscillator}.
In this case the condition of non existence of ''continuous flow'' is that  $3P_k-(2\times 3+1)(P_{k-1})\ge 0
 \Leftrightarrow P_{k-1}\le \frac{3P_k}{7}$.
Hence repeating the same process if  $P_{i-1}\le \frac{3P_i}{7}$ ($0\le i\le k$), 
there will always remain a blank interval and the ''flow of 1'' would not 
be continuous.
Since for all $0\le i\le k$ one has    $P_i\le \lfloor\frac{l_i}{5}\rfloor +4$, it follows that the condition    
   $\lfloor\frac{l_{i-1}}{5}\rfloor \hskip -.1 cm +4\le \frac{3}{7}$ $(\lfloor\frac{l_i}{5}$ $\rfloor +2)$ 
(iq4) also implies 
  that the ''flow of 1'' is not continuous. 
Remark that we take a lower bound for $P_i$ and upper bound for $P_{i-1}$ which leads to a stronger condition 
on the minimum size $l_m$. 
Then by simplification of (iq4) we get that $3l_i-7l_{i-1}\ge 110$.
Using (iq1)$\equiv (l_{i-1}\le \lfloor\frac{3l_i}{10}\rfloor +6$) that gives the minimum condition for a train generated by an
 {\it oscillator} $\mathcal{C}_{l_i}$ to cross the following one $\mathcal{C}_{l_{i-1}}$
 we obtain  that if $3l_i-7(\lfloor\frac{3l_i}{10}\rfloor +6)\ge 110$ (iq5) then the condition 
 $3l_i-7l_{i-1}\ge 110$ (iq4) remains true.
 The simplification of inequality (iq5) leads to $l_{i}\ge \frac{1520}{9}\approx 168.88$.
Since we have chosen $l_{\bf m}=170$ as the minimum size of the {\it oscillators}, 
there is no equicontinuous points in $(S(\mu_c),F_e)$  which finish the proof.

 \hfill$\Box$

The dynamical system ($X$,$F_e$) has  equicontinuous points  since patterns of the type $E^*0^kE_0$ 
(with $E^*\in\hat{E}$)  will 
produce a continuous flow of letters 1.
Note that we can construct  $F'_e$ a CA similar to $F_e$ with a more complex local rule such that 
the initial measure $\mu'_I$ is the uniform measure on $X$ and the invariant measure obtained by Ces\`{a}ro means 
is similar to the measure $\mu_c$ used in our example.
   In this case $(X,F'_e,\mu_c )$ and $(S(\mu_c),F'_e,\mu_c)$ both have 
   $\mu_c$-equicontinuous points but no equicontinuous  points.

\subsection{Questions}

-Is it possible to find a sensitive, $\mu$-equicontinuous and $\mu$-invariant CA when $\mu$ is 
the uniform measure? Or more generally when the topological support $S(\mu )$ 
of the  $F$-invariant measure is a mixing subshift of finite type?




-To simplify the proof of Proposition \ref{noneq}, we have taken a lot of upper bounds (
 for example the use of (iq1) is very strong).
We wonder   what is the minimum size  for the {\it oscillators} 
 such that $(S (\mu ),F_e)$   can not produce equicontinuous points?
 
-From inequality (iq3) in Proposition \ref{pro10}, it is possible to see that for each  
point $x$ and  $m\in\NN$ such that $\mu_c(B_m(x))>0$,  we have $B_m^{\mu_c}(x)=B_m(x)$.
What are the conditions on the minimum size of the counters $l_m$ in order to loose this property?
 In this case is it possible that there is no equicontinuous points? 
  Remark that for some points $x$, $B_m^{\mu}(x)\neq B_m(x)$ in the Gilman's example of $\mu$-equicontinuous 
CA given in the beginning of section 3.3.

 -Is there exists a $\mu$-invariant and $\mu$-equicontinuous CA  such that 
  $(S(\mu ),F)$ has no equicontinuous points and  there exists $m\in\NN$ and a point  $x\in A^\ZZ$ with   
  $B_m(x)\neq B_m^\mu (x)$ and $\mu (B_m(x))>0$?

-More generally, what type of dynamic characterize sensitive and $\mu$-equi\-con\-ti\-nuous CA $(S(\mu),F)$  
and how common is this behavior that seems to appear in different
simulations of one dimensional CA?

\section*{Acknowledgments}
 
The author whish to acknowledge the NSERC Discovery Grant \#562620, the CNPq and the Department of Mathematics at Trent University 
in which a part of the work have been done.

\end{document}